\documentclass[11pt]{article}

\usepackage{amsmath}
\usepackage{amssymb}
\usepackage{eucal}

\begin{document}

\baselineskip 16pt

\title{Finite groups with permutable  Hall  subgroups}

\author{Xia Yin, Nanying Yang\thanks{Research  is supported by
a NNSF grant of China (Grant \#11301227) and Natural Science Foundation of Jiangsu Province (grant \# BK20130119).}\\
{\small School of Science, Jiangnan University}\\ {\small Wuxi 214122 , P. R. China}\\
{\small E-mail:yangny@jiangnan.edu.cn}\\ \\
}

 \date{}
\maketitle

\begin{abstract}  Let $\sigma =\{\sigma_{i} | i\in I\}$ be a partition of
 the set of all primes $\Bbb{P}$ and  $G$  a finite group.
A set  ${\cal H}$  of subgroups of $G$ is said to be a  \emph{complete Hall
 $\sigma $-set} of $G$   if   every  member $\ne 1$ of
 ${\cal H}$ is a Hall $\sigma _{i}$-subgroup  of $G$ for some $i\in I$ and
 $\cal H$ contains exactly one Hall  $\sigma _{i}$-subgroup of $G$
 for every  $i$ such that $\sigma _{i}\cap   \pi (G)\ne \emptyset$.

 In this paper, we study the
structure of $G$ assuming that some subgroups of $G$ permutes
with all members of ${\cal H}$.

\end{abstract}


\let\thefootnote\thefootnoteorig

\section{Introduction}

Throughout this paper, all groups are finite and $G$ always denotes
a finite group. We use $\pi (G)$ to denote the set of all primes
dividing $|G|$.
A subgroup $A$ of $G$ is said to \emph{permute} with a subgroup $B$
if $AB=BA$.  In this case they  say also that the subgroups $A$ and $B$ are
\emph{permutable}.

Following \cite{1}, we  use $\sigma$  to denote some partition of
$\Bbb{P}$. Thus   $\sigma =\{\sigma_{i} |
 i\in I \}$, where   $\Bbb{P}=\cup_{i\in I} \sigma_{i}$
 and $\sigma_{i}\cap
\sigma_{j}= \emptyset  $ for all $i\ne j$.

A  set  ${\cal H}$ of subgroups of $G$ is a
 \emph{complete Hall $\sigma $-set} of $G$ \cite{commun, PI}  if
 every member $\ne 1$ of  ${\cal H}$ is a Hall $\sigma _{i}$-subgroup of $G$
 for some $\sigma _{i} \in \sigma$ and ${\cal H}$ contains exactly one Hall
 $\sigma _{i}$-subgroup of $G$ for every  $i$ such that $\sigma _{i}\cap  \pi(G)\ne \emptyset$.
  If    every two members of   ${\cal H}$ are permutable,
 then ${\cal H}$  is said to be a \emph{$\sigma$-basis} \cite{2}
of $G$. In the case when   ${\cal H} =\{\{2\}, \{3\},  \ldots \}$
 a complete Hall $\sigma $-set  ${\cal H}$ of $G$
 is also called  \emph{a complete set  of  Sylow subgroups} of $G$.

We use  $\frak{H}_{\sigma}$ to denote the class of all soluble
groups $G$ such that   every  complete Hall $\sigma $-set of $G$
 forms a $\sigma$-basis of $G$.

A large number of publications are connected with study the
situation when   some subgroups of $G$ permute  with all members of
some fixed complete set of Sylow subgroups of $G$. For example, the
 classical  Hall's result  states: {\sl G is soluble if and only if it has a
Sylow basis, that is, a complete set of pairwise permutable Sylow
subgroups}. In \cite{h2} (see also Paragraph 3 in \cite[VI]{hupp}),
Huppert proved that $G$ is a soluble group in which  every complete set of
Sylow subgroups forms a Sylow basis if and only if the automorphism
group induced by $G$ on every its chief factor $H/K$
 has the  order divisible by at most one different from $p$
prime, where $p\in \pi (H/K)$. In the paper  \cite{h1},  Huppert  proved that if  $G$ is
soluble and it has a complete set $\cal S$ of
 Sylow subgroups such that every maximal subgroup of every subgroup
in $\cal S$ permutes with all other members of $\cal S$, then $G$ is
supersoluble.

The above-mentioned  results in \cite{h2, hupp, h1} and many other related
results make natural to ask:

(I) {\sl Suppose that $G$ has a complete Hall $\sigma $-set
$\cal H$  such that  every maximal subgroup of any subgroup in
${\cal H} $ permutes with all other members of ${\cal H}$.  What we
can say then about the structure of $G$? In particular, does it
true then that $G$ is supersoluble in the case when every member of ${\cal H}$ is supersoluble?}

(II) {\sl  Suppose that $G$ possesses a a complete Hall $\sigma $-set.
 What we can say then about the structure of $G$ provided
every   complete Hall $\sigma $-set of $G$ forms  a $\sigma$-basis in
$G$?}

Our first observation is the following result concerning Question (I).

{\bf Theorem A.} {\sl Suppose that $G$ possesses a a complete Hall
$\sigma $-set $\cal H$ all  whose members are supsersoluble. If   every maximal
subgroup of any non-cyclic subgroup in ${\cal H} $ permutes with all
other members of ${\cal H}$, then   $G$ is supersoluble.}

In the classical case, when $\sigma =\{\{2\}, \{3\}, \ldots \}$, we get
from Theorem A the following  two known results.

 {\bf Corollary 1.1} (Asaad M., Heliel  \cite{asaad}).  {\sl If  $G$
 has a complete set $\cal S$ of
 Sylow subgroups such that every maximal subgroup of every subgroup
in $\cal S$ permutes with all other members of $\cal S$, then $G$ is
supersoluble. }

Note that  Corollary 1.1 is proved in \cite{asaad} on the base of the
classification of all simple non-abelian groups. The proof of Theorem A
does not use such a classification.

 {\bf Corollary 1.2} (Huppert  \cite[VI, Theorem 10.3]{hupp}).
 {\sl If every Sylow subgroup of $G$ is cyclic,
 then $G$ is  supersoluble. }

The class $1\in \mathfrak F$ of groups  is said to be a \emph{ formation} provided
every homomorphic image of $G/G^{\mathfrak F}$ belongs to $\mathfrak
F$. The formation $\mathfrak F$ is said to be: \emph{saturated}
provided $G\in \mathfrak F$ whenever $G^{\mathfrak F}\leq \Phi (G)$;
\emph{hereditary} provided $G\in \mathfrak F$ whenever $G\leq A\in
{\mathfrak F}$.

Now let  $p > q > r$ be primes such that $qr$ divides $p-1$. Let $P$
be a group of order $p$ and  $QR\leq Aut (P)$, where $Q$ and $R$ are
groups with order $q$ and $r$, respectively. Let $G=P\rtimes (QR)$.
Then, in view of the above-mentioned Hupper's result in \cite{h2},
$G$ is not a  group such that every complete set of Sylow subgroups
forms a Sylow basis of $G$. But it is easy to see that every complete Hall
 $\sigma$-set of $G$, where $\sigma= \{\{2, 3 \},  \{7 \},  \{2, 3, 7\}'\}$,
 is a $\sigma$-basis of $G$. This
elementary example is a motivation for our next result, which gives
the answer to Question (II) in the universe    of all soluble
groups.

{\bf Theorem B.} {\sl  The class $\frak{H}_{\sigma}$ is a hereditary
formation and it  is saturated if and only if $|\sigma| \leq 2$.
Moreover,  $G\in {\frak{H}}_{\sigma}$ if and only if $G$ is soluble
and the automorphism group induced by $G$ on every its chief factor
of order divisible by $p$  is either a $\sigma_{i}$-group, where
$p\not\in \sigma _{I}$, or a  $(\sigma _{i}\cup \sigma _{j})$-group for
some different  $\sigma _{i} $  and $ \sigma _{j}$ such that $p\in \sigma _{i}$.}

In the  case when $\sigma =\{\{2\}, \{3\}, \ldots \}$ we get
from Theorem B the following

 {\bf Corollary 1.3} (Huppert \cite{h2} ).  {\sl Every complete set of
Sylow subgroups of a soluble group $G$ forms a Sylow basis of $G$ if and only if the automorphism
group induced by $G$ on every its chief factor $H/K$
  has   order divisible by at most one different from $p$
prime, where $p\in \pi (H/K)$. }

\section{Proof of Theorem A}

{\bf Lemma 2.1} (See  Knyagina and  Monakhov \cite{knyag}). {\sl
Let $H$, $K$  and $N$ be pairwise permutable
subgroups of $G$ and  $H$ is a Hall subgroup of $G$. Then $N\cap HK=(N\cap H)(N\cap K).$}

{\bf Proof of Theorem A.}  Assume that this theorem is false and let $G$ be a counterexample of
minimal order.  Let  ${\cal H}=\{H_{1}, \ldots , H_{t} \}$.
 We can assume, without loss of generality, that  the smallest
prime divisor $p$ of $|G|$ belongs to $\pi (H_{1})$.  Let $P$ be a
Sylow $p$-subgroup of $H_{1}$.

(1) {\sl If  $R$ is a minimal normal subgroup of $G$, then $G/R$  is
supersoluble. Hence  $R$ is the   unique minimal normal subgroup of $G$, $R$
 is not cyclic and $R\nleq \Phi (G)$}.

We show that the hypothesis holds for $G/R$.  First note that
$${\cal H}_{0}=\{H_{1}R/R, \ldots , H_{t}R/R \}$$ is
  a  complete Hall
$\sigma$-set    of $G/R$, where
$H_{i}R/R\simeq H_{i}/H_{i}\cap R$ is supersoluble since $H_{i}$ is supersoluble by hypothesis
for all $i=1, \ldots , t$.

 Now let $V/R$ be a maximal subgroup
of $H_{i}R/R$, so $|(H_{i}R/R):(V/R)| =p$ is a prime.
   Then    $V=R(V\cap H_{i})$ and hence
 $$p=|(H_{i}R/R):(V/R)|=  |(H_{i}R/R):(R(V\cap H_{i})/R)|=|H_{i}R:R(V\cap
H_{i})|=$$ $$=
|H_{i}||R||R\cap (V\cap H_{i})|:|V\cap
H_{i}||R||H_{i}\cap R|=
|H_{i}|:|V\cap H_{i}|=|H_{i}:(V\cap H_{i})|,$$ so  $V\cap H_{i}$ is
 a maximal subgroup of $H_{i}.$  Assume that  $H_{i}R/R$ is not cyclic.
Then $H_{i}$ is not cyclic, so $$(V\cap H_{i})H_{j}=H_{j}(V\cap H_{i})$$ for all $j\ne i$
    by hypothesis and hence   $$(V/R)(H_{j}R/R)=(R(V\cap H_{i})/R)(H_{j}R/R)=
(H_{j}R/R)((V\cap H_{i})R/R)=(H_{j}R/R)(V/R).$$
  Consequently the hypothesis holds for $G/R$,
so  $G/R$ is supersoluble by  the choice of $G$. Moreover, it is well
known that the class of all supersoluble groups is a saturated formation
(see Ch. VI in \cite{hupp} or ??? in \cite{Guo-bookI}). Hence the choice of $G$ implies that
$R$ is the unique minimal normal subgroup of $G$, $R$
 is not cyclic and $R\nleq \Phi (G)$.

(2) {\sl $G$ is not soluble. Hence $R$ is not abelian and $2\in \pi
(R)$.}

Assume that this is false. Then $R$ is an abelian $q$-group for some
prime $q$. Let $q\in \pi_{k}$. Since $R$ is non-cyclic  by Claim (1) and $R\leq H_k$,
$H_{k}$ is  non-cyclic. Hence every member of
${\cal H}$ permutes with each maximal subgroup of $H_k$. Since
$R\nleq \Phi (G)$, $R\nleq \Phi (H_{k})$ and so there exists a
maximal subgroup $V$ of $H_{k}$ such that $R\nleq  V$ and $RV=H_{k}$.
Hence $E=R\cap V\ne 1 $ since $|R| > q$ and  $H_k$ is supersoluble.
Clearly, $E$ is normal in $H_k$. Now assume that $i\ne k$. Then $V$
permutes with $H_{i}$ by hypothesis, so $VH_{i}$ is a subgroup of
$G$  and   $$R\cap VH_{i}=(R\cap V)(R\cap H_{i})=R\cap V=E$$ by Lemma 2.1
 and so $H_i\leq N_G(E).$
Therefore $H_{i} \leq  N_{G}(E)$ for all $i=1, \ldots , t.$ This
implies that $E$ is normal in $G$, which contradicts the minimality
of $R$. Hence we have (2).

(3) {\sl If $R$  has a Hall $\{2, q\}$-subgroup for each $q$
dividing $|R|$, then a  Sylow 2-subgroup $R_2$ of $R$ is
non-abelian.}

Assume that this is false. Then by Claim (2) and Theorem 13.7 in
\cite[XI]{BlH}, the composition factors of $R$ are isomorphic to one
of the following groups: a) $PSL(2, 2^{f})$;  b) $PSL(2, q)$, where
8 divides $q-3$ or $q-5$; c) The Janko group $J_{1}$; d) A Ree
group. But with respect to  each of these groups it is well-known
(see, for example \cite[Theorem 1]{Tyut}) that the group has no a
Hall $\{2, q\}$-subgroup for at least one odd prime $q$ dividing its
order. Hence we have (3)

(4) {\sl If at least one of the  subgroups  $H_{i}$ or  $H_{k}$, say $H_{i}$, is
non-cyclic,  then  $H_{i}H_{k}= H_{k}H_{i}$ } (This
follows from the fact that every maximal subgroup of $H_{i}$
permutes with $H_{k}$).

(5)  {\sl $H=H_{1}$ is not cyclic} (This directly follows from Claim (2),
 \cite[IV, 2.8]{hupp} and the Feit-Thompson theorem).

In view of Claim (5), ${\cal H}$   contains  non-cyclic subgroups.
Without loss of generality, we may assume that $H_{1}, \ldots ,
H_{r}$ are non-cyclic groups  and all  groups $H_{r+1}, \ldots ,
H_{t}$ are cyclic.

(6) {\sl Let  $E_{\{i, j\}}=H_{i}H_{j}$ where $i \leq r$. If  $r$ is
the smallest prime dividing $|E_{\{i, j\}}|$,  then $E_{\{i, j\}}$
is $p$-nilpotent, so it is soluble. Therefore $E_{\{i, j\}}\ne G$.}

Clearly, the hypothesis holds for $E_{\{i, j\}}$. Hence if $E_{\{i,
j\}} < G$, then this subgroup is supersoluble by the choice of $G$,
and so it is $p$-nilpotent. Now assume that $E_{\{i, j\}}=G$. Then
$r=p=2$ and $E_{\{i, j\}}=HH_{j}=H_jH$. Let $V_{1}, \ldots , V_{t}$
be the set of all maximal subgroups of a Sylow $2$-subgroup $P$ of
$H$. Since  $H$ is supersoluble,  it has a normal $2$-complement
$S$. Then $SV_{i}$ is a maximal subgroup of $H$, so
$SV_{i}H_{j}=H_{j}SV_{i}$ is a subgroup of $G$ by hypothesis.
Moreover, this subgroup is normal in $G=E_{\{i, j\}}$ since
$|G:H_{j}SV_{i}|=2$. Now let $E=SV_{1}H_{j} \cap \cdots \cap
SV_{t}H_{j}$. Then $E$ is normal in $G$  and clearly $E\cap P\leq
\Phi (P)$.

Now we show that for any prime $q$ dividing $|H_{j}|$, there
are a Sylow $q$-subgroup $Q$ of $H_{j}$ and an element $h\in H$ such
that $P\leq N_{G}(Q^{h})$. Indeed, by the Frattini argument,
$G=EN_{G}(Q)$. Hence by \cite[VI, 4.7]{hupp}, there are Sylow
$2$-subgroups $G_2$, $E_2$ and $N_2$ of $G$, $E$ and $N_{G}(Q)$
respectively such that $G_2=E_2N_2$. Let $P=(G_2)^{x}$. Then
$P=(E_{2})^{x}(N_{2})^{x}$, where $(E_2)^{x}$ is a Sylow
$2$-subgroup of $E$ and $(E_{2})^{x}$ is a Sylow $2$-subgroup of
$(N_{G}(Q))^{x}=N_{G}(Q^{x})$. Since $G=HH_{j}$, $x=hw$ for some
$h\in H$ and $w\in H_{j}$. Hence
$$N_{G}(Q^{x})=N_{G}(Q^{wh})=N_{G}((Q^{w})^{h}),$$ where $Q^{w}$ is a
Sylow $q$-subgroup of $H_{j}$. Therefore $(E_{2})^{x}=E\cap P\leq
\Phi (P)$. Consequently, $P\leq N_{G}((Q^{w})^{h})$. This shows that
for any prime $q$ dividing $|H_{j}|$, there is a Sylow $q$-subgroup
$Q$ of $H_{j}$ and an element $h\in H$ such that $P\leq
N_{G}(Q^{h})$. Thus $G$  has a Hall $\{2, q\}$-subgroup $PQ^{h}$ for
each $q$ dividing $|H_{j}|$. Moreover, since $H$ is supersoluble by
hypothesis, $G$ has a Hall $\{2, s\}$-subgroup for each $s$ dividing
$|H|$. Hence in view of Claim (3), $P$ is not abelian. Then $P\cap
F(H)\ne 1$, so $P\cap F(H)\leq Z_{\infty}(H)$ since $H$ is
supersoluble. Let $Z$ be a group of order $2$ in $Z(H)$. Since
$Z\leq P\leq N_{G}((Q^{h})$, $Z=Z^{h^{-1}}\leq N_{G}(Q)$. It follows
that $Z\leq N_{G}(H_{j})$. Thus $Z^{G}=Z^{HH_{j}}=Z^{H_{j}}\leq
ZH_{j}$. This shows that a Sylow $2$-subgroup of $Z^{G}$ has order
$2$. Hence $Z^{G}$ is $2$-nilpotent. Let $S$ be the $2$-complement
of $Z^G$. It is clear that $S\ne 1$. Since $S$ is characteristic in
$Z^G$, it is normal in $G$. On the other hand, $S$ is soluble by the
Feit-Thompson theorem. This induces that $G$ has an abelian minimal
normal subgroup, which contradicts Claim (2). Thus (6) holds.

(7) {\sl    $E_{i}=HH_{i}$ is supersoluble  for all $i =2, \ldots ,
t$} ((Since the  hypothesis holds for  $E_{i}$ and $E_{i} < G$ by
Claim (5),  this follows from the choice of $G$).

(8) {\sl $E= H_{1} \ldots H_{r}$ is soluble.}

We argue by induction on $r$. If $r=2$, it is true by Claim (5). Now
let  $r > 2$ and assume that the assertion is true for $r-1$. Then
by Claim (4), $E$ has at least three   soluble  subgroups $E_{1}$,
$E_{2}$, $E_{3}$ whose indices  $E:E_{1}|$, $|E:E_{2}|$, $|E:E_{3}|$
are pairwise coprime. But then $E$ is soluble by the Wielandt
theorem \cite[I, 3.4]{DH}.

(9) {\sl $R$  has a Hall $\{2, q\}$-subgroup for each $q$ dividing
$|R|$.}

It is clear in the case when $q\in \pi (H)$. Now assume that $q\in
\pi (H_{i})$ for some $i > 1$. Then Claim (6) implies  that
$B=HH_{i}$ is a Hall soluble subgroup of $G$. Hence $B$ has a Hall
$\{2, q\}$-subgroup $V$ and so $V\cap R$ is a Hall $\{2,
q\}$-subgroup of $R$.

(10) {\sl A Sylow 2-subgroup $R_2$ of $R$ is non-abelian} (This
follows from Claims (3) and (9)).

(11) {\sl If $q\in \pi (H_{k})$ for some $k > r$, then $q$ does not
divide $|R:N_{R}((R_2)')|$.}

By Claim (7), $B=HH_{k}$ is supersoluble. Hence  there is a Sylow
$q$-subgroup of $Q$ of $B$ such that $PQ$ is a Hall $\{2, q
\}$-subgroup of $B$. Then $U=PQ\cap R=(P\cap R)(Q\cap R)=R_2(Q\cap
R)$ is a Hall supersoluble subgroup of $R$ with cyclic Sylow
$q$-subgroup $Q\cap R$. By \cite[VI, 9.1]{hupp}, $Q\cap R$ is normal
in $U$, and $U/C_{U}(Q\cap R)$ is an abelian group by \cite[Ch. 5,
4.1]{Gor}. Hence
$$R_2C_{U}(Q\cap R)/C_{U}(Q\cap R)\simeq R_2/R_2\cap C_{U}(Q\cap R)$$ is
abelian and so $(R_2)'\leq C_{U}(Q\cap R)$. Consequently, $Q\cap R\leq
N_{R}((R_2)')$.

{\sl The final contradiction.}
In view of Claim (11), $R=(E\cap R)   N_{R}((R_2)')$. Hence
$$((R_2)')^{R}=((R_2)')^{(E\cap R)   N_{R}((R_2)')}=((R_2)')^{E\cap
R}\leq E\cap R.$$  But by Claim (8), $E\cap  R$ is soluble. On the other
hand, Claim (10) implies that $(R_2)'\ne 1$ and  so $R$ is soluble,
contrary to Claim (2).
The theorem is thus proved.

\section{Proof of Theorem B}

The following lemma can be proved by the direct calculations on the
base of well-known properties of Hall subgroups of soluble
subgroups.

 {\bf Lemma 3.1.} {\sl The class
$\frak{H}_{\sigma}$ is closed under taking homomorphic images,
subgroups and direct products}.

{\bf  Proof of Theorem B.} Firstly, from Lemma 3.1, $\frak{H}_{\sigma}$ is a
hereditary formation.

Now we prove that $G\in {\frak{H}}_{\sigma}$ if and only if $G$ is
soluble and the automorphism group induced by $G$ on every its chief
 factor of order divisible by $p$  is either a $\sigma_{i}$-group, where
$p\not\in \sigma_{i}$, or a  $(\sigma_{i}\cup
\sigma_{j})$-group for some different $\sigma_{i}$ and $ \sigma_{j}$
 such that    $p \in \sigma_{i}$.

{\sl Necessity.} Assume that this is false and let $G$ be a counterexample of
 minimal order. Then $G$ has a chief factor $H/K$ of order divisible by $p$ such that
$A=G/C_{G}(H/K)$ is neither a $\sigma_{i}$-group, where
$p\not\in \sigma_{i}$, nor a $(\sigma_{i}\cup
 \sigma_{j})$-group, where $\sigma_{i}\ne  \sigma_{j}$ and  $p \in \sigma_{i}$.
 Since $$G/C_{G}(H/K)\simeq
(G/K)/(C_{G}(H/K)/K) =(G/K)/C_{G/K}(H/K)
$$ and the hypothesis hods for  $G/K$ by Lemma 3.1, the choice of $G$ implies that    $K=1$.

First we show that  $H\ne C_{G}(H)$. Indeed, assume that
$H=C_{G}(H)$. By  hypothesis, every complete Hall  $\sigma $-set  ${\cal
W}=\{W_{1}, \ldots , W_{t} \}$ of $G$    forms a
$\sigma $-basis of  $G$. Without loss of generality, we can
assume that $p\in \pi(W_{1})$. It is cleat that  $ t > 2$. Since
$H=C_{G}(H)$,  $H$ is the unique minimal normal subgroup of $G$ and
$H\nleq \Phi (G)$ by \cite[Ch.A, 9.3(c)]{DH} since $G$ is soluble.
 Hence $H=O_{p}(G)=F(G)$ by
\cite[Ch.A, 15.6]{DH}. Then for some maximal subgroup $M$ of $G$ we
have $G=H\rtimes M$. Let $V=W_{3}$.  We now show that $V^{x}\leq
C_{G}(W_{2})$ for all  $x\in G$. First note that
$W_{2}V^{x}=V^{x}W_{2}$ is a Hall $(\sigma_{2}\cup \sigma_{3})$-subgroup of
$G$. Since $|G:M|$ is a power of $p$,  any Hall $\sigma_{0}$-subgroup
of $M$,  where $p\not \in \pi_{0}$, is a Hall $\pi_{0}$-subgroup of
$G$. Hence we can assume without loss of generality that
$W_{2}V^{x}\leq M$ since $G$ is soluble.  By hypothesis,
$W_{2}(V^{x})^{y}=(V^{x})^{y}W_{2}$ for all $y\in G$, so $$D=\langle
(W_{2})^{V^{x}} \rangle \cap \langle (V^{x})^{W_{2}}\rangle$$ is
subnormal in $G$ by \cite[1.1.9(2)]{prod}. But $D\leq \langle W_{2},
V^{x} \rangle \leq M$,
 so  $$D^{G}=D^{HM}=D^{M}\leq M_{G}=1$$
 by \cite[Ch. A, 14.3]{DH}, which implies that $[W_{2}, V^{x}]=1$.
  Thus $V^{x}\leq
C_{G}(W_{2})$ for all  $x\in G$. It follows that
 $H\leq (W_{3})^{G}\leq N_{G}(W_{2})$ and therefore $W_{2}\leq C_G(H)=H$, a
 contradiction. Hence $H\ne C_{G}(H)$.

Finally,   let $D=G\times G$, $A^{*}=\{ (g, g)|g\in G \}$, $C=\{ (c,
c)|c\in C_{G}(H) \}$ and $R=\{(h,1)|h\in H \}$.  Then $C\leq
C_{D}(R)$,  $R$ is a minimal normal subgroup of $A^{*}R$ and the factors $R/1$
 and $RC/C$ are $(A^{*}R)$-isomorphic.
 Moreover,
 $$C_{A^{*}R}(R)=R(C_{A^{*}R}(R)\cap A^{*})=RC,$$ so
$$A^{*}R/C=(RC/C)\rtimes (A^{*}/C),$$
 where $A^{*}/C\simeq A$ and
$RC/C$ a minimal normal subgroup of $A^{*}R/C$ such that
$C_{A^{*}R/C}(RC/C)=RC/C$. As $H < C_G(H)$, we see that
$|A^{*}R/C| < |G|$. On the other hand, by Lemma 3.1, the hypothesis
holds for $A^{*}R/C$, so  the choice of $G$ implies that $A\simeq A^{*}/C$
  is  either a $\sigma_{i}$-group, where
$p\not\in \sigma_{i}$, or   a $(\sigma_{i}\cup
 \sigma_{j})$-group for some different $\sigma_{i}$ and $ \sigma_{j}$
 such that    $p \in \sigma_{i}$.
This contradiction completes the proof of the  necessity.

{\sl Sufficiency.}   Assume that this is false and let $G$ be a
counterexample of minimal order.  Then $G$ has a  complete Hall set
 ${\cal W}=\{W_{1}, \ldots  , W_{t} \}$
 of type $\sigma$ such that for some  $i$ and $j$
 we have  $W_{i} W_{j}\ne W_{j} W_{i} $. Let $R$ be a minimal normal subgroup of
$G$. Then:

(1) {\sl    $G/R\in {\frak{H}}_{\sigma}$, so $R$ is a unique minimal
normal subgroup of $G$.}

It is clear that the hypothesis holds for $G/R$, so $G/R\in
{\frak{H}}_{\sigma}$ by the choice of $G$. If $G$ has a minimal
normal subgroup $L\ne R$, then we also have $G/L\in
{\frak{H}}_{\sigma}$. Hence $G$ is isomorphic to some subgroup of
$(G/R)\times (G/L)$ by \cite[I, 9.7]{hupp}. It follows from Lemma 3.1 that $G\in
{\frak{H}}_{\sigma}$. This contradiction shows that we
have Claim (1).

(2) {\sl The hypothesis holds for any subgroup $E$ of $G$.}

Let $H/K$ be any chief factor of $G$ of order divisible by $p$ such
that $H\cap E\ne K\cap E$. Then $G/C_{G}(H/K)$ is  either a $\sigma_{i}$-group, where
$p\not\in \sigma_{i}$, or   a  $(\sigma_{i}\cup
\sigma_{j})$-group   for some different  $\sigma_{i}$ and $\sigma_{j}$ such that
  $p \in \sigma_{i}$. Let $H_{1}/K_{1}$ be a chief factor
of $E$ such that $K\cap E \leq K_{1} < H_{1} \leq H\cap E$. Then
$H_{1}/K_{1}$ is a $p$-group and
$$EC_{G}(H/K)/C_{G}(H/K)\simeq E/ (E\cap C_{G}(H/K))$$ is  either a $\sigma_{i}$-group or  a
$(\sigma_{i}\cup  \sigma_{j})$-group.
Since  $$C_{G}(H/K)\cap
E\leq C_{E}(H\cap E/K\cap E)\leq C_{E}(H_{1}/K_{1}),$$
$E/C_{E}(H_{1}/K_{1})$ is also  either a $\sigma_{i}$-group or  a $(\sigma_{i}\cup
\sigma_{j})$-group.  Therefore the hypothesis holds for
every factor $H_{1}/K_{1}$ of some  chief series of $E$. Now
applying the Jordan-H\"older Theorem for chief series we get Claim (2).

(3) {\sl    $R$ is   a Sylow $p$-subgroup of $G$.}

Since $G/R\in  {\frak{H}}_{\sigma}$ by Claim (1),   $$(W_{i}R/R)
(W_{j}R/R)=(W_{j}R/R)(W_{i}R/R),$$ so  $W_{i} W_{j}R$ is a subgroup
of $G$. Assume that $R$ is not a Sylow $p$-subgroup of $G$ and let
$B=W_{i} W_{j}R$. Then $B\ne G$. On the other hand, the hypothesis
holds for $B$ by Claim (2). The choice of $G$ implies that
$B\in {\frak{H}}_{\sigma}$, so $W_{i} W_{j}= W_{j} W_{i} $, a
contradiction. Hence Claim (3) holds.

{\sl Final contradiction for sufficiency.}
In view of Claims  (1) and (3),  there is a
maximal subgroup $M$ of $G$ such that $G=R\rtimes M$ and $M_{G}=1$.
Hence $R=C_{G}(R)=O_{p}(G)$ by \cite[Ch.A, 15.6]{DH}.  Since $p$
does not divide $|G:R|=|G:C_{G}(R)|$ by Claim (3),  the hypothesis implies that
 $M\simeq G/R$ is  a Hall
 $\sigma_{k}$-group
 for some  $ \sigma_{k}\in \sigma $, so one of the subgroups $W_{i}$
 or $ W_{j}$ coincides with $R$. Thus
 $G=W_{i} W_{j}= W_{j} W_{i}. $
This contradiction completes the proof of the sufficiency.

Finally we prove that $\frak{H}_{\sigma}$ is saturated if and only if $|\sigma |\leq 2.$
It is clear
that  $\frak{H}_{\sigma}$ is a saturated  formation for any $\sigma$
with $|\sigma|\leq 2$. Now we  show that for any $\sigma$ such
that $|\sigma|> 2$, the formation $\frak{H}_{\sigma}$ is not
saturated.

Indeed, since $|\sigma| > 2$, there are primes  $p < q < r$ such
that for some distinct $\sigma _{i}$, $\sigma _{j}$ and $\sigma _{k}$ in
$\sigma$ we have  $p\in \sigma _{i}$, $q\in \sigma _{j}$ and $r\in \sigma
_{k}$. Let $C_{q}$ and $C_{r}$ be groups of order $q$ and $r$,
respectively. Let $P_{1}$  be a simple $\Bbb{F}_{p}C_{q}$-module
which is faithful for $C_{q}$, $P_{2}$  be a simple
$\Bbb{F}_{p}C_{r}$-module which is faithful for $C_{r}$. Let
$H=P_{1}\rtimes C_{q}$ and $Q$ be a simple $\Bbb{F}_{q}H$-module
which is faithful for $H$. Let $E=(Q\rtimes H)\times (P_{2}\rtimes
C_{r})$.

Let $A=A_p(E) $ be the $p$-Frattini module of $E$
(\cite[p.853]{DH}), and let $G$ be a non-splitting extension of $A$
by $E$. In this case, $A\subseteq \Phi (G)$ and $G/A\simeq E.$ Then $G/\Phi (G)\in
 \frak{H}_{\sigma}$, where $\sigma=\{\sigma_i, \sigma_j, \sigma _k\}.$
By Corollary 1 in \cite{{Griess}}, $QP_{1}P_{2}=O_{p',
p}(E)=C_{E}(A/Rad(A))$. Hence for some normal subgroup $N$ of $G$ we
have $A/N\leq \Phi (G/N)$ and $G/C_{G}(A/N)\simeq C_{q}\times
C_{r}$ is a $(\sigma _i\cup \sigma_j)$-group. But neither
 $p\not\in \sigma_i$ nor $p\in \sigma_j$.
 Hence $G\not\in \frak{H}_{\sigma}$ by
the necessity.
The theorem is proved.

\end{document}